# A SUPPLEMENT TO THE BOSE-DASGUPTA-RUBIN (2002) REVIEW OF INFINITELY DIVISIBLE LAWS AND PROCESSES

*By* S SATHEESH

*Telecom Training Centre*, *Trichur*.

ssatheesh@sancharnet.in

*SUMMARY*. This supplement stresses that if a discrete distribution $\{p_i,\ i = 0,1,2,\ \ldots\}$ on $I_0$ is infinitely divisible with integer-valued components then $p_0>0$, and discusses some of its implications. We then give certain recent developments and references not reported in the Bose-Dasgupta-Rubin (2002) review in *Sankhya-A*, and some new results in the topics (i) infinitely divisibility and stability of discrete laws, (ii) random infinite divisibility, (iii) operator stable laws, (iv) class-L laws, (v) Goldie-Steutel result, (vi) max-infinite divisibility and stability, (vii) simulation, (viii) alternate stable laws, (ix) applications and (x) free probability theory.

## 1. Introduction

The Bose-Dasgupta-Rubin (2002) (abbreviated as BDR (2002) here) paper extensively reviews the literature on infinitely divisible (ID) laws and processes with many illustrative examples and a list of Levy measures. This supplement is motivated by the need to (i) stress that when a discrete distribution $\{p_i,\ i = 0,1,2,\ \ldots\}$ on $I_0 = \{0,\ 1,\ 2,\ \ldots\}$ is ID with integer-valued components, $p_0>0$ is necessary (ii) explain why ID laws on $I_0$ cannot have any gaps (iii) discuss a class of discrete laws useful in transfer theorems to study schemes with random sample size and (iv) present certain recent developments and references not reported in BDR (2002). Topics mentioned in the summary will be discussed in the following sections. New results/ ideas stressed is given as Theorem/ Lemma/ Remark/ Note/ Example.

## 2. Discrete laws that are infinitely divisible, stable and in class-L

In discussions on ID r.vs $X$ on $I_0$ with integer-valued components $\{X_i\}$, one comes across the phrase "an ID lattice r.v $X$ with $P\{X{=}0\} > 0$ ". See eg. Katti (1967), Steutel (1973, 1979). So a natural question is whether there could be an ID r.v $X$ for which $P\{X{=}0\} = 0$. Feller (1968, p.290) while proving that ID laws on $I_0$ must be compound Poisson (and conversely) states this condition as an analytical requirement in terms of the corresponding probability generating function (PGF). Steutel and van Harn (1979) assume this condition in the statement of the result itself. Bondesson (1981) while reviewing ID laws discusses those on $I_0$ and also those obtained by truncating out the mass at zero. But whether they continue to be ID or not is not mentioned. Johnson *et al.* (1992, p.352) proves

---

that if $X$ is compound Poisson (that is, ID) and has finite mean then P$\{X=0\}$>0 and a partial converse of this. BDR (2002) also does not discuss this point.

Hence first we give a probabilistic argument to prove that the condition P$\{X=0\}$>0 is necessary for all ID laws on $I_0$ with integer-valued components. (I am thankful to Professor M Sreehari, M S University of Baroda, India, for suggesting the line of argument presented here. Also the author has come to know that Kallenberg has priority to this result, see Grandall (1997, p.26)). Among other implications it explains why ID laws on $I_0$ cannot have any gaps; c.f Remark.9 of BDR (2002).

**Theorem 1** If a r.v $X$ on $I_0$ is ID with integer-valued components then necessarily P$\{X = 0\} > 0$.

**Proof.** Under the assumptions on $X$, there exists i.i.d r.vs $\{X_i\}$ on $I_0$ such that

$$X \stackrel{d}{=} X_1 + \ldots + X_n \quad \text{for every } n \geq 1. \tag{2.1}$$

Now assume the contrary that P$\{X=0\} = 0$. Then P$\{X_i=0\} = 0$ for all $1 \leq i \leq n$ and every $n \geq 1$. Let $k$>0 be the least integer such that P$\{X = k\}$>0. For a given $n \geq 1$, let $r$ >0 is the least integer such that P$\{X_i = r\}$>0. Then the minimum value assumed with positive probability by the RHS of (2.1) is $nr$ where as that of the LHS is $k$ and they are never equal when $n$>$k$. Hence when $n$>$k$ we cannot have the representation (2.1). The minimum values $nr$ and $k$ are equal for every $n \geq 1$ only if $k = 0$ = $r$. But $k = 0$ implies $r = 0$. Hence P$\{X = 0\} > 0$. Also P$\{X_i=0\} > 0$ for all $i = 1, \ldots, n$.

One can now see that the discussion on shifted zero truncated Poisson law in Kemp (1978,p.1436) is not at all required to show that it is not ID.

**Theorem 2** If a r.v $X$ on $I_0$ is ID, then P$\{X$ >$k\}$> 0 for any positive integer $k$. In other words P$\{X$ >$k\} = 0$ for some positive integer $k$ is impossible.

The proof follows by a similar line of argument. These two theorems imply the following:

(i) The support of an ID law on $I_0$ cannot be a finite subset of $I_0$. Recall the well-known property that the support of ID laws cannot be finite. It is now clear that the binomial law is not ID.

(ii) From Theorem.1, if $Q(s)$ is a PGF that is ID, then $Q(0) \neq 0$. This can be seen, as the discrete analogue of the fact the characteristic function (CF) of ID laws does not have real zeroes.

For an ID r.v $X$ on $I_0$, P$\{X = 0\} = 0$ implies that the components of $X$ are no more supported by positive integers. Thus for any PGF $Q(s)$, the PGF $P(s) = sQ(s)$ is not ID with integer-valued components. Eg. The geometric law on $I_1 = \{1, 2, \ldots\}$ is not ID with integer-valued components though the geometric law on $I_0$ is ID with integer-valued components. But we know that the property of infinite divisibility is invariant under translation. Hence the geometric law on $I_1$ is also ID but its components are not integer-valued. Now, can we distinguish the class of ID laws on $I_0$ whose components are integer-valued and those whose components are not? Dr Iksanov has suggested the following sharpening of theorem.1. See also Kallenberg's result referred to by Grandall (1997, p.26).

**Theorem 1*a*** Let $X$ be a non-negative integer-valued ID r.v. Then $P\{X=0\}>0$ iff the support of $X$ coincides with that of its components for every $n \geq 1$ integer.

**Proof.** Theorem.1 proves the if part and lemma.2.2 of Hu *et al.* (2002) the only if part.

**Remark 1** With reference to the Remark.9 of BDR (2002) one may ask; Why there are no gaps in the support of a discrete ID r.v $X$ with $P\{X=1\}>0$? For the ID r.v $X$ (with integer-valued components) we now have both $P\{X=0\}>0$ and $P\{X=1\}>0$. Hence for every $n \geq 1$ the components $X_i$ also satisfy $P\{X_i=0\}>0$ and $P\{X_i=1\}>0$. Consequently $X$ cannot have any gaps in its support as the summation (2.1) ensures that every non-negative integer carries a probability mass.

An interesting question now is; are there any discrete ID r.vs with integer-valued components, that have gaps in its support? Of course, they cannot have a mass at $X=1$. The answer is in the affirmative and we have:

**Example 1*a*** From the PGF of the negative binomial law we have the PGF $\{p/(1-qs^k)\}^t$, $p+q=1$, $k>1$ integer and $t>0$. This is the PGF of the negative binomial sum of r.vs degenerate at $k>1$ integer. Obviously this PGF is ID and its probability masses are $k-1$ integers apart.

**Example 1*b*** Again let $X$ be a r.v with PGF $s\{p/(1-qs^k)\}^t$, $p+q=1$, $k>1$ integer and $t>0$. We can see that: here $X$ is ID, $P\{X=0\}=0$, $P\{X=1\}>0$ and $X$ has gaps in its support, but the components are not integer-valued.

**Remark 1*a*** Thus remark.9 of BDR (2002) holds good only for ID laws on $I_0$ with $P\{X=0\}>0$.

The components of an ID law in the sense of (2.1) are themselves ID and hence their support is unbounded. On the other hand recall that any ID law is generated by an 'infinite number' of uniformly small components. In this sense the components need not be ID and hence be boundedly supported or the support of the component variables could be finite.

Eg. A Poisson law can be written as a compound of Bernoulli laws as:

$$\exp\{-\lambda(1-s)\} = \exp\{-ab(1-s)\} = \exp\{-a[1-(1-b+bs)]\}, \text{ where } ab = \lambda \text{ and } 0<b<1.$$

Now we derive the result of Feller (1968, p.290) that ID laws on $I_0$ must be compound Poisson (and conversely). We do this in two parts for their independent interest.

**Theorem 3** If an ID law on $[0,\infty)$ has an atom at the origin then it must be compound Poisson.

**Proof.** Let $\varphi(s) = e^{-\psi(s)}$ be the Laplace transform (LT) of the ID law (Feller, 1971, p.450). If a distribution with LT $\varphi$ has an atom at the origin then it is equivalent to $\varphi(\infty)>0$. For an ID law this is reflected as $\psi(\infty) = \lambda < \infty$. Setting $\psi(s)/\lambda = F(s)$, we have $F(o) = 0$ and $F(\infty) = 1$. Further since $\psi(s)$ has completely monotone derivative (CMD), $F(s)$ is continuous and non-decreasing. Thus $F(s)$ is a d.f with CMD. Hence by Pillai and Sandhya (1990) $F(s)$ is the d.f of a mixture of exponential laws and hence $F(s) = 1 - \alpha(s)$, where $\alpha(s)$ is a LT. Hence $\varphi(s) = \exp\{-\lambda[1-\alpha(s)]\}$ completing the proof.

**Corollary 1** A distribution on $[0,\infty)$ with LT $\varphi$ is compound Poisson iff $-\log\varphi$ is a d.f.

**Theorem 4** Discrete analogue of ID laws on $[0,\infty)$ are Poisson compounds of non-negative integer valued r.vs.

**Proof.** Readily follows from Theorems.1 & 3 and noticing that the components must be discrete.

We have come to know that a different proof of this result is available in Ospina and Gerber (1987) and its multivariate extension by Sundt (2000).

Steutel and van Harn (1979) after describing the discrete class-L ($L$) and discrete stable laws by their PGFs, recorded that consideration of a formal discrete analogue of $L$ laws led them to compound Poisson, ie. ID laws. Using a method for constructing distributions on the non-negative lattice of points $I_o = \{0,1,2, ....\}$ as discrete analogue of continuous distributions on $[0,\infty)$, Satheesh (2001a), (see also Satheesh and Nair (2002a) and Satheesh *et al.* (2002)) have justified the notions of discrete $L$ and discrete stable laws by deriving these PGFs from the LT of their continuous counter parts. They also developed discrete analogue of distributions of the same type (D-type) and discussed the role of Bernoulli law in the context. Their key ideas are:

**Lemma 1** If $\varphi(s)$ is a LT, then $Q(s) = \varphi(1-s)$, $0<s<1$ is a PGF. Conversely, if $Q(s)$ is a PGF and $Q(1-s)$ is completely monotone for all $s>0$, then $\varphi(s) = Q(1-s)$ is a LT.

**Definition 1** Two PGFs $Q_1(s)$ and $Q_2(s)$ are of the same D-type if $Q_1(1-s) = Q_2(1-cs)$, for all $0<s<1$, or equivalently, $Q_1(u) = Q_2(1-c+cu)$ for all $0<u<1$, and some $0<c<1$.

**Theorem 5** Two non-negative lattice r.vs $X$ and $Y$ will have the same D-type distribution if and only if

$$X = \sum_{i=1}^{Y} Z_i \quad \text{for some i.i.d Bernoulli r.vs } \{Z_i\} \text{ independent of } Y, \text{ the Bernoulli probability being } c.$$

This theorem also justify the replacement of $cX$ in the continuous set up by $c \cdot X = \sum_{i=1}^{Y} Z_i$ to obtain the corresponding lattice analogue, as done in Steutel and van Harn (1979). Subsequently, from the definition of $L$ laws the Steutel-van Harn definition of discrete $L$ laws follows by invoking Lemma.1 and definition.1. Discrete stable laws with PGF $\exp\{-\lambda(1-s)^\alpha\}$ is readily obtained from the LT $\exp\{-\lambda(s)^\alpha\}$ of $\alpha$-stable laws. Satheesh and Nair (2002a) and Satheesh *et al.* (2002) discusses the summation stability properties of discrete semi stable laws in general. Another explanation why the normalizing parameter for stable (continuous & discrete) should be $n^{1/\alpha}$, $\alpha \in (0.1]$ is also given in Satheesh *et al.* (2002). Some other works in this area are: Rao and Shanbhag (1994, p.160), Christoph and Schreiber (1998), Anil (2001) and Sajikumar (2002).

Now, why should stable laws be absolutely continuous? (Feller, 1971, p.215). What precisely is the requirement for a distribution to be stable? By Feller (1971, p.176) stable distributions are ID and distinguished by the fact that $X_i$ differs from $X$ only by location parameters. One is tempted to ask these questions because the sum of any number of Poisson laws is again Poisson. What is the ground for the argument that the Poisson law obtained as the sum is different from the component Poisson by more than location parameters or they do not satisfy the normalizing requirement (Romano and Seigel, 1986, Ex.4.41, p.81) ?. Of course, by location parameters, Feller means distributions of the same type (Feller, 1971, p.45, 137, 169, 176). Here we try to answer these two questions, though it is well accepted that stable laws are absolutely continuous. Notice also the similarity between the PGFs of discrete stable laws and the LTs of $\alpha$-stable laws on [0,∞).

Stable family is a location-scale family of laws such that if random variables $X$ and $Y$ each have distributions from that family then the distribution of their sum, i.e. $X+Y$, also belongs to the same family. More precisely one may call this as summation stability under linear normalization. Notice that in the definition of stable laws the requirement is: $X$ is stable if for each positive integer $n$ there exists real constants $b_n > 0$ and $a_n$ such that with $\{X_i\}$ i.i.d as $X$ and $S_n = X_1 + …. + X_n$,

$$S_n \overset{d}{=} b_n X + a_n.$$

One can very well see that such a normalizing requirement would never be satisfied by a discrete law as the support immediately changes. So the culprit is the definition of types that is only applicable to continuous laws. This is a case of over precision given by the mathematical description of an idea. In the discrete case the nature of linear normalization is given by definition.1. Poisson laws and in general distributions with PGF $\exp\{-\lambda(1-s)^\alpha\}$, $\alpha \in (0.1]$ are stable. In place of Poisson law as the counter example discussed in Romano and Seigel (1986, p.81) on can consider the Laplace law, which is not stable but is ID (in fact is in $L$).

A characteristic property of stable laws is that a distribution possesses a domain of attraction iff it is stable. For distributions on [0,∞) in terms of LTs, this is equivalent to (Feller,1971, p.448):
$-n \log \varphi(s/a_n) \longrightarrow -\log \varsigma(s)$. Now invoking our Lemma.1 we do have PGFs $\varphi((1-s)/a_n)$ and $\varsigma(1-s)$, and they satisfy $-n \log \varphi((1-s)/a_n) \longrightarrow -\log \varsigma(1-s)$. Accordingly, $\varsigma(1-s) = \exp\{-(1-s)^\alpha\}$, $0 < \alpha \leq 1$ which is discrete stable. Hence:

**Theorem 6** In the discrete case, discrete stable laws and they alone possess domain of attraction.

Thus stable laws are not necessarily absolutely continuous. It is the manifestation or nature of linear normalization that is different in the continuous and discrete supports.

More on ID laws are available also in: Kallenberg (1975), Riedel (1980*a*, *b*) and Seshadri (1993, chapter.5).

### 3. Random Infinite Divisibility and Stability

From the angle of random summation, generalizations of ID laws were first considered by: Klebanov, *et al.* (1984) for geometric sums, introducing geometrically ID (GID) laws. The book by Kalashnikov (1997) and two recent reviews by Kozubowski and Rachev (1999*a,b*) may be consulted for more on geometric summation schemes. Other works not reported therein are: Sandhya (1991*a*), see also Sandhya and Pillai (1999) - introducing and studying attraction and partial attraction in geometric sums, Pillai (1990) – studying harmonic mixtures and GID laws, Pillai and Sandhya (1990) – GID laws and mixtures of exponentials, Sandhya (1991*b*) – GID laws, Cox processes and p-thinning of renewal processes, generalizing the Renyi characterization of Poisson processes and showing that a Cox process is renewal iff its interval distribution is GID, Pillai and Sandhya (1996) - geometric sums in reliability. Levy measures in the context of GID laws were discussed in Pillai and Sandhya (2001) characterizing inverse Gaussian law. A discussion of GID laws and integrated Cauchy functional equation is available in Pillai and Anil (1997). Geometric summation schemes, tailed distributions and time series models have been considered by Pillai (1991), Pillai, *et al.* (1995), Pillai and Jayakumar (1994), Jayakumar and Pillai (1993), Jose and Pillai (1996), Jose and Seethalakshmi (1999) and Balakrishna and Jayakumar (1996, 1997), Seethalekshmi and Jose (2001, 2002, 2003, 2004), Jose and Alice Thomas (2002, 2003), Alice Thomas and Jose (2003*a,b*). In the discrete case geometric sum stability was defined and discussed in Jayakumar (1995), and in more generality by Satheesh (2001*a*) and Satheesh and Nair (2002*a*), see also Pillai and Jayakumar (1995).

Satheesh (2001a) and Satheesh *et al.* (2002) have discussed generalization of stability of geometric sums by studying distributions that are stable under summation w.r.t Harris law. They showed that the notion of stability of random sums can be extended to include the case when $X$ is discrete and proposed a method to identify the probability law of $N$ for which $X$ is N-sum stable. The following result can be easily proved.

**Theorem 7** Absolutely continuous laws cannot be represented as an N-sum if P$\{N=0\}>0$. In particular they can never be N-sum stable w.r.t an $N$ with P$\{N=0\} > 0$.

Satheesh *et al.* (2002) showed that in the discrete case this is possible and gave an example.

Further generalization of ID and GID laws to random sums were considered in the following works: In the N-ID laws in Sandhya (1991a, 1996), the description was not constructive but she discussed two examples of non-geometric laws for $N$. The description of N-ID laws (with CF $\varphi(\psi)$, where exp$\{-\psi\}$ is the CF of an ID law) by Gnedenko and Korolev (1996), Klebanov and Rachev (1996) and Bunge (1996), are based on the assumption that the PGF $\{P_\theta, \theta \in \Theta\}$ of $N$ formed a commutative semi-group. Satheesh (2001*b*, 2002) showed that this assumption is not natural in a generalization of ID laws since it captures the notion of stable rather than ID laws. Also, it rules out any $P_\theta$ having an atom at the origin. Two glaring situations are that this theory cannot (i) handle compound Poisson law which is at the heart of the concept of ID laws and (ii) approximate negative binomial sums by taking $\varphi$ as gamma as one expects. For exact N-sums Satheesh, *et al.* (2002) showed that when $\varphi$ is gamma, $N_0$ is Harris($a,k$) and Satheesh (2001$a$) showed that the limit law of Harris as $\theta \to 0$ is indeed N-ID with $\varphi$ being gamma. Kozubowski and Panorska (1996, 1998) discussed $\nu$-stable laws with CF $\varphi\{\psi\}$ (exp{-$\psi$} being a stable CF) approximating $N_0$-sums, assuming $N_0 \xrightarrow{\ p\ } \infty$ as $\theta \downarrow 0$ and it could handle Poisson and negative binomial sums. But how to identify those $N_0 \xrightarrow{\ p\ } \infty$ ?

The following two lemmas are ramifications of Feller's proof of Bernstein's theorem (Feller, 1971, p.440) and were first proved by Satheesh (2001$b$) see also Satheesh (2002). Together they provide a class of discrete laws that can be used in the transfer theorems for sums and extremes Gnedenko (1982). They also suggest why the property of the classical sum in (2.1) is reflected as a Poisson sum. Because the limit law of sums as $n$ increases has the CF exp{-$\psi$} (that of an ID law) and the limit law of $\theta N_0$, as $\theta \to 0$, where $N_0$ is Poisson($\theta$) has LT $e^{-s}$.

**Lemma 2** $\wp_\varphi(s) = \{P_0(s) = s^j \varphi\{(1-s^k)/\theta\}$, $0<s<1$, $j \geq 0$ & $k \geq 1$ integer and $\theta > 0\}$ describes a class of PGFs for any given LT $\varphi$.

**Lemma 3** Given a r.v $U$ with LT $\varphi$, the integer valued r.vs $N_0$ with PGF $P_\theta$ in the class $\wp_\varphi(s)$ described in Lemma.2 satisfy $\theta N_0 \xrightarrow{\ d\ } kU$ as $\theta \to 0$.

Notice that the class $\wp_\varphi(s)$ (of Lemma.2) cater for a wide variety of discrete laws with support on $I_0$ or $I_1$ or any countable subclass of it permitting gaps as well. Satheesh (2001$b$), motivated by the limitations in the description of random infinite divisibility already stated, introduced $\varphi$-ID laws with CF $\varphi\{\psi\}$ approximating $N_0$-sums where the PGF of $N_0$ is a member of $\wp_\varphi(s)$ without requiring $N_0 \xrightarrow{\ p\ } \infty$. This generalizes the fact that the limit distributions of compound Poisson laws are ID. See also Satheesh (2002) from the angle mixtures of ID laws and Satheesh and Sandhya (2002).

**Definition 2** A CF $f(t)$ is $\varphi$-ID if for every $\theta \in \Theta$, there exists a CF $g_\theta(t)$, a PGF $P_\theta$ that is independent of $h_\theta$, such that $P_\theta\{g_\theta(t)\} \to f(t) = \varphi\{-\log \omega(t)\} \; \forall \; t \in \mathbf{R}$, as $\theta \downarrow 0$ through a $\{\theta_n\} \in \Theta$. Here, $\varphi$ is the LT of a r.v $Z>0$ and $\omega(t)$ is a CF that is ID.

For the class of PGFs $\wp_\varphi(s)$: (i) Satheesh (2001$b$) proved a result, that is analogous to the Theorem 4.6.5 in Gnedenko and Korolev (1996, p.149), and discussed $\varphi$-attraction and partial $\varphi$-attraction and (ii) Satheesh and Sandhya (2002) showed that a N&S condition for the convergence of $N_0$-sums of $\{X_{\theta,j}\}$ with CF $g_\theta$ to a $\varphi$-ID law is $(1 - g_\theta(t))/\theta \to \psi(t)$ as $\theta \downarrow 0$ through a $\{\theta_n\}$, where $\psi(t)$ is a continuous function. This result generalizes Theorem.1.1 of Feller (1971, p.555). Analogous result for random vectors is also given therein. From these two results one can conceive the notions of attraction and partial attraction for $N_0$-sums when the PGF of $N_0$ belongs to $\wp_\varphi(s)$. Here we will denote the PGFs by $P_n$ -corresponding to $\theta_n = 1/n$ , $\{n\}$ the sequence of positive integers and $P_{n_m}$ - corresponding to $\theta_m = 1/n_m$ , $\{n_m\}$ a subsequence of $\{n\}$. See also Satheesh (2002).

**Definition 3** A CF $g(t)$ belongs to the domain of φ-attraction (Dφ-A) of the CF $f(t)$ if there exist sequences of real constants $a_n = a(\theta_n) > 0$ and $b_n = b(\theta_n)$ such that with $g_n(t) = g(t/a_n)$ exp($-itb_n$); $\underset{n\to\infty}{Lt}\ P_n\{g_n(t)\} = f(t)$, $\forall t \in R$, and $g(t)$ belongs to the domain of partial φ-attraction (DPφ-A) of $f(t)$ if; $\underset{m\to\infty}{Lt}\ P_{n_m}\{g_m(t)\} = f(t)$, $\forall t \in R$.

Certain implications of the two results are: They enable us to conclude that if the CF $g(t)$ belongs to the DA of the stable law (DPA of the ID law) with CF $\omega(t) = e^{-\psi(t)}$, then it is also a member of Dφ-A of the φ-stable law (DPφ-A of a φ-ID law) with CF $f(t) = \varphi\{\psi(t)\}$ and the converses are also true. All that we need is to prescribe, $\theta_n = 1/n$ for φ-attraction and $\theta_m = 1/n_m$ for partial φ-attraction. Thus the DA of a stable law (DPA of an ID law) with CF $\omega(t) = e^{-\psi(t)}$ coincides with the Dφ-A of the φ-stable law (DPφ-A of a φ-ID law) with CF $f(t) = \varphi\{\psi(t)\}$ for each PGF $P_\theta \in \mathscr{P}_\varphi(s)$, and none of them are empty as well. These notions and conclusions generalize those on attraction and partial attraction in geometric sums in Sandhya (1991), Sandhya and Pillai (1999), and those for N-sums in Gnedenko and Korolev (1996) and Klebanov and Rachev (1996). One may extend these results to random vectors on $\boldsymbol{R}^d$ for $d \geq 2$ integers. This extends and generalizes proposition 2.1 of Kozubowski and Panorska (1998). See also Satheesh (2002).

**Remark 2** The flexibility achieved by our φ-ID laws over the N-ID and/ or ν-stable laws is that we have a class of distributions for the random sample size for every LT $\varphi$. Using the structure of compound Poisson laws one can extend the notion of N-ID, ν-stable and φ-ID laws to the discrete setup.

**Remark 3** In N-ID or φ-ID laws, a curiosity is whether their infinite divisibility has anything to do with that of $N$ or $\varphi$. Satheesh (2001$b$) showed with examples, that (i) an N-ID law can be ID even when $N$ is not ID (ii) a φ-ID law can be ID even when $\varphi$ is not ID.

## 4. Operator Stable Laws

Some other works on operator stable laws are Meerschacrt and Scheffler (2001). Let $X$, $X_1$, $X_2$, ….are i.i.d random vectors on $\boldsymbol{R}^d$ with a common distribution $\mu$ and $Y_0$ is a random vector on $\boldsymbol{R}^d$ whose distribution $\tau$ is full, that is, not supported on any lower dimensional hyperplane. Then $\tau$ is operator stable if there exists linear operators $A_n$ on $\boldsymbol{R}^d$ and a sequence of points $b_n \in \boldsymbol{R}^d$ such that

$$A_n \sum_{i=1}^{n}(X_i - b_n) \xrightarrow{\ d\ } Y_0 \ \text{ or } \ A_n\mu_n * \varepsilon(s_n) \xrightarrow{\ w\ } \tau , \ \varepsilon(s_n) \text{ being the mass at the point } -nA_nb_n.$$

In the above setup let $Y$ be a random vector on $\boldsymbol{R}^d$ whose distribution $\lambda$ is full. Let $Z > 0$ be a r.v with probability distribution $\nu$ and suppose there exists integer valued r.vs $\{N_n\}$ independent of $\{X_i\}$ such that $N_n/n \xrightarrow{d} Z$. Then $\lambda$ or $Y$ is operator $\nu$-stable (Kozubowski, *et al.* (2002)) if there exists linear operators $A_n$ on $\boldsymbol{R}^d$ and a sequence of points $b_n \in \boldsymbol{R}^d$ such that

$$A_n \sum_{i=1}^{N_n} (X_i - b_n) \xrightarrow{d} Y \text{ or } f(t) = \int_0^\infty \omega(t)^s \, d\nu(s) \, ; \, f(t) \text{ \& } \omega(t) \text{ being the CFs of } Y \text{ and } Y_0 \, .$$

When $N_n$ are geometric(1) with mean $n$, $Z$ is standard exponential then $\lambda$ or $Y$ is operator geometric-stable with $f(t) = 1/(1 - \log\{\omega(t)\})$.

In view of Lemmata 2 & 3 we extend the notion to operator $\varphi$-stable laws as follows. Let $Y$ be a random vector on $\boldsymbol{R}^d$ with CF $f(t)$ whose distribution $\lambda$ is full. Let $Z > 0$ be a r.v with probability distribution $\nu$ and LT $\varphi$. Then:

**Definition 4** The CF $f(t)$ (or the distribution $\lambda$ or the random vector $Y$) is operator $\varphi$-stable if there exists linear operators $A_\theta$ on $\boldsymbol{R}^d$ and a sequence of points $b_\theta \in \boldsymbol{R}^d$ such that as $\theta \downarrow 0$ through $\{\theta_n\}$

$$A_\theta \sum_{i=1}^{N_\theta} (X_i - b_\theta) \xrightarrow{d} Y \text{ or } f(t) = \int_0^\infty \omega(t)^s \, d\nu(s) = \varphi\{-\log \omega(t)\}.$$

Let $\{X_{\theta,j} \, ; \, \theta \in \Theta, j \geq 1\}$ be a sequence of i.i.d random vectors on $\boldsymbol{R}^d$ and for $k \geq 1$ integer, set $S_{\theta,k} = X_{\theta,1} + \ldots + X_{\theta,k}$. Now we will consider a $\{\theta_n \in \Theta\}$ and as $\theta \downarrow 0$ through $\{\theta_n\}$ $S_{\theta,n} \xrightarrow{d} U$ so that $U$ is operator stable. Also in an $N_\theta$-sum of $\{X_{\theta,k}\}$, $N_\theta$ and $X_{\theta,k}$ are assumed to be independent for each $\theta \in \Theta$. From Lemmata.2 & 3 and the transfer theorem we then have:

**Theorem 8** The limit law of $N_\theta$-sums of $\{X_{\theta,j}\}$ as $\theta \downarrow 0$ through a $\{\theta_n\}$, where the PGF of $N_\theta$ is a member of $\mathscr{P}_\varphi(s)$, is necessarily operator $\varphi$-stable. Conversely, for any given LT $\varphi$; the operator $\varphi$-stable law can be obtained as the limit law of $N_\theta$-sums of i.i.d random vectors as $\theta \downarrow 0$ for each member of $\mathscr{P}_\varphi(s)$.

**Remark 4** We can now have definitions of generalized domains of $\varphi$-attraction of operator $\varphi$-stable laws and results analogous to the $\varphi$-stable laws generalizing those in Kozubowski, *et al.* (2002). The flexibility achieved by our operator $\varphi$-stable laws over the operator $\nu$-stable laws is mentioned in the applications in remark.6.

## 5. Class-L Laws

The following references are worth consulting in the context of class-L ($L$) laws. These laws are important as they can be used in the construction of ARMA processes, Gaver and Lewis (1980). There are many methods to verify whether a distribution is a member of $L$ or not, see eg. Lukacs (1970), Shanbhag and Sreehari (1977, 1979), Shanbhag *et al.* (1977), Ismail and Kelkar (1979), Pillai and Satheesh (1992), Pillai and Sabu George (1984), Sandhya and Satheesh (1996*a*) and another review with examples and Levy measures in Jurek (1997).

Mixtures of infinitely divisible (ID) laws were introduced by Feller (1971, p.573) and studied by Keilson and Steutel (1974) and under the name φ-mixtures by Satheesh (2002), $\varphi$ being the LT of the mixing distribution. Sandhya and Satheesh (1996), Ramachandran (1997) and Satheesh, *et al.* (2002) have discussed the relation between random sum stability and $L$ laws. From the angle of φ-mixtures of ID laws Satheesh (2002) generalize Theorem.2.1 (part (iii), (iv) and (v)) of Ramachandran (1997) and Theorems 2.3, 2.4 and 2.5 regarding generalized Linnik laws of Satheesh, *et al.* (2002) and Erdogan and Ostrowski (1998). Construction of $L$ laws was considered in Sandhya and Satheesh (1996*a*), Satheesh *et al.* (2002) and Satheesh (2001, 2002*a*). The result is:

**Theorem 9** A φ-mixture of a strictly stable law is in $L$ if $\varphi \in L$.

For more on $L$ laws and their relation to random processes see: Barndorff-Neilson and Shephard (2000), Diedhiou (1998), Iksanov and Jurek (2003), Jeanblanc *et al.* (2001), Jurek (2000, 2002, 2003), Knight (2001) and Sato (2001). Generalizations of discrete-$L$ laws are considered by Dimitrov and Kolev (2002). Pillai and Jayakumar (1994) has discussed a generalized class-L property.

## 6. Goldie-Steutel Result

Though ID laws have found fruitful utility in both theory and applications it was a bit hard to verify whether or not a given distribution is ID. Goldie (1967) achieved a breakthrough in this direction for non-negative continuous r.vs. He showed that distributions that are mixtures of exponential laws are ID. Steutel (1969) went further by showing that densities of mixtures of exponential laws are completely monotone (CM) and hence all CM densities are ID. This result brought in the notion of CM functions to the realm of probability densities. Shanbhag and Sreehari (1977, 1979) and Shanbhag, et. al (1977) extended Goldie's result to include mixtures of gamma and discussed their relation to $L$ laws. By the method of Goldie, Thorin (1977) observed that Pareto laws are ID and that they are in $L$ also thus bringing in the notion of generalized gamma convolutions that are in $L$. Pillai and Sandhya (1990) strengthened the Goldie-Steutel result by showing that CM densities are characteristic of mixtures of exponentials and further they are GID. Thorin had proved this result independently in the context of stationary renewal processes in Insurance mathematics assuming the mixing distribution to have finite mean, where they have found a lot of applications, see Grandall (1991, p.49). Pillai and Sandhya (1990) also gave a density that is GID but not CM. Some more results on mixtures of exponential laws are available in Sandhya and Satheesh (1996*b*, 1997). Bondesson (1990) discusses recent results on generalized gamma convolutions and CM functions. The discrete analogue of mixtures of exponentials as mixtures of geometric laws on $I_0$ has been discussed by Satheesh and Sandhya (1997). Sandhya and Satheesh (1997) showed that a mixed Poisson process is renewal iff it is Poisson improving on McFadden (1965) who proved it assuming the mixing law to have finite mean.

## 7. Max-Infinite Divisibility and Stability

Stability of the maximum and/ or minimum with random (*N*) sample size was discussed by Voorn (1987,1989) who introduced the concept for continuous laws. Pillai (1991) and Pillai and Sandhya (1996) have discussed the notion with samples of geometric size and characterized semi-Pareto laws. See also Shaked (1975) and Sreehari (1995). Discussing the important case of the exponential law, Satheesh (2001*a*) (see also Satheesh and Nair (2002*b*)), has further extended the notion to include discrete laws as well. We just give the definition to discuss a counter example.

**Definition 5** Two lattice laws with d.fs $F(k)$ and $G(k)$ are of the same type if $G(k) = F(\alpha k)$ for all $k = 0$, 1, …. , and for some $\alpha > 0$.

In the discrete case this is possible as follows. Let $F(k) = P\{X < k\} = 1 - m(k)$, $k = 0, 1, ….$ and $G(k) = P\{Y < k\} = 1 - m(\alpha k)$, $k = 0, 1, ….$ where $\{m(k)\}$ is the sequence of realizations of a LT $m(s)$, $s>0$, at the non-negative integral values of $s$. These d.fs can also be viewed from the angle of mixtures of geometric laws on $I_0$ and in this case $\{m(k)\}$ is the moment sequence of the mixing law .

Notice that Definition.5 appears different from Definition.1 but is quite similar to the one in the continuous case. In the continuous case the ideas are equivalent. But :

**Example 2** Let *X* has a geometric(0,*p*) law.

Then its d.f is $F(k) = 1 - q^k$, $k = 0,1,2, ….$ , $q = 1 - p$ and PGF is $Q_X(s) = p/(1 - qs)$.

Now in accordance with Definition.5 consider the r.v *Y* with d.f $G(k) = 1 - q^{ck}$ for some $0 < c < 1$.

Setting $q = \frac{1}{4}$ and $c = \frac{1}{2}$, we have $q^c = \sqrt{(\frac{1}{4})} = \frac{1}{2}$. Further;

$Q_X(s) = 3/(4 - s)$ and $Q_Y(s) = 1/(2 - s)$ and

$Q_X(1 - \frac{1}{2} + s/2) = 6/(7 - s)$ and $Q_Y(\frac{1}{2} + s/2) = 2/(3 - s)$.

Thus neither $Q_X(s) = Q_Y(\frac{1}{2} + s/2)$ nor $Q_Y(s) = Q_X(\frac{1}{2} + s/2)$ considering both the possibilities. Hence the two definitions are not equivalent.

Sreehari (1995) and Satheesh (2001*a*) (also Satheesh and Nair (2002*b*)) proposed a method to identify the probability law of $N$ for which $X$ is N-max/ N-min stable. A conjecture characterizing the geometric law on $I_1$ in this context is given by Satheesh (2001*a*), also Satheesh and Nair (2000), explaining certain observations made by Arnold, *et al.* (1986) and Marshall and Olkin (1997).

Again these are stability under linear normalization. Under power normalization this has been touched upon in Sreehari (1995).

Paralleling the classical ID and stable laws Balkema and Resnick (1977) discussed max-infinite divisibility and stability (MID and max-stable laws) for d.fs on $\boldsymbol{R}^d$ for $d \geq 2$ integer. An asymptotic distribution of extremes with random sample size has been studied by Silvestrov and Teugels (1998) and Dorea and Goncalves (1999). Rachev and Resnick (1991) and Mohan (1998) have extended the notion of MID laws to the geometric sample size paralleling GID laws. Satheesh (2002) and Satheesh and Sandhya (2002) have further extended this to random sample sizes with PGFs in $\mathscr{P}_\varphi(s)$ paralleling φ-ID laws, and discussed partial φ-max-attraction for φ-MID laws with d.f $\varphi\{-\log H\}$, where $H$ is MID, generalizing the results for the geometric case given in Rachev and Resnick (1991) and Mohan (1998). Let $\{\boldsymbol{Y}_{\theta,j} ; \theta \in \Theta, j \geq 1\}$ be a sequence of i.i.d random vectors in $\boldsymbol{R}^d$, $d \geq 2$ integer with d.f $H_\theta$ and for $k \geq 1$ integer, set $\boldsymbol{M}_{\theta,k} = \vee \{\boldsymbol{Y}_{\theta,1}, \ldots., \boldsymbol{Y}_{\theta,k}\}$. Again we will consider a $\{\theta_n \in \Theta\}$ and $\{k_n ; n \geq 1\}$ of natural numbers such that as $\theta \downarrow 0$ through $\{\theta_n\}$ and $\boldsymbol{M}_{\theta,k_n} \xrightarrow{d} \boldsymbol{V}$ so that $\boldsymbol{V}$ is MID. Also in an $N_\theta$-max of $\{\boldsymbol{Y}_{\theta,j}\}$, $N_\theta$ and $\boldsymbol{Y}_{\theta,j}$ are independent for each $\theta \in \Theta$.

**Definition 6** A d.f $F$ is φ-MID if for every $\theta \in \Theta$, there exists a d.f $H_\theta$, a PGF $P_\theta$ that is independent of $H_\theta$, such that $P_\theta\{H_\theta\} \to F = \varphi\{-\log G\}$ as $\theta \downarrow 0$ through a $\{\theta_n\}$. Here, $\varphi$ is the LT of a r.v $Z > 0$ and $G$ is a d.f that is MID.

**Theorem 10** The limit law of $N_\theta$-maxs of i.i.d random vectors as $\theta \downarrow 0$, where the PGF of $N_\theta$ is a member of $\mathscr{P}_\varphi(s)$, is necessarily φ-MID. Conversely, for any given LT $\varphi$; the φ-MID law is the limit law of $N_\theta$-maxs of i.i.d random vectors as $\theta \downarrow 0$ for each member of $\mathscr{P}_\varphi(s)$.

Satheesh and Sandhya (2002) showed that a N&S condition for the convergence of $N_\theta$-maxs of $\{\boldsymbol{Y}_{\theta,j}$ with d.f $H_\theta\}$ to a φ-MID law is $(1 - H_\theta(t))/\theta \to \mu$ as $\theta \downarrow 0$ through a $\{\theta_n\}$, where $\mu$ is an exponent measure. This motivated them the following definition and subsequent conclusions.

**Definition 7** A d.f $H$ in $\boldsymbol{R}^d$, $d \geq 2$ integer, belongs to the Dφ-MA of the d.f $F$ with non-degenerate marginal distributions if there exists normalizing constants $a_{i,n} = a_i(\theta_n) > 0$ and $b_{i,n} = b_i(\theta_n)$ such that with $H_n(\boldsymbol{y}) = H(a_{i,n}y_i + b_{i,n}, 1 \leq i \leq d)$; $\underset{n \to \infty}{Lt} P_n\{H_n\} = F$ and $H$ belongs to the DPφ-MA of the d.f $F$ if;

$\underset{m \to \infty}{Lt} P_{n_m}\{H_m\} = F.$

Thus if the d.f $H$ belongs to the DMA of a max-stable law (DPMA of a MID law) with d.f $G = e^{-\mu}$ then it is also a member of the D$\varphi$-MA of a $\varphi$-max-stable law (DP$\varphi$-MA of a $\varphi$-MID law) with d.f $F = \varphi\{\mu\}$ and the converses are also true. All that we need is to prescribe, $\theta_n = 1/n$ for $\varphi$-max-attraction and $\theta_m = 1/n_m$ for partial $\varphi$-max-attraction. Thus the DMA of a max-stable law (DPMA of a MID law) with d.f $G = e^{-\mu}$ coincides with the D$\varphi$-MA of a $\varphi$-max-stable law (DP$\varphi$-MA of a $\varphi$-MID law) with d.f $F = \varphi\{\mu\}$ for each PGF $P_0 \in \mathscr{P}_\varphi(s)$, and none of them are empty. See also Satheesh (2002).

**Remark 5** The flexibility achieved by our $\varphi$-MID laws over other discussions of similar limit laws is that we have a class of distributions for the random sample size for every LT $\varphi$.

## 8. Simulation

Algorithms for generating non-negative integer valued r.vs from PGFs or a moment sequence has been discussed in Devroye (2001). Stochastic representation of a r.v in terms of other r.vs is the key to some other methods of simulation. Kozubowski (2000) and Kozubowski and Panorska (1999) make use of this approach to generate univariate and multivariate geometric stable and limit laws in random-sum schemes. Generation of $\varphi$-stable laws can be done in this approach, as the stochastic representation for these laws is the same as those for the corresponding $\nu$-stable laws.

## 9. Alternate Stable Laws

As we have already discussed the notion of stable laws is precisely the property of invariance of a class of location-scale family of laws under the operation of summation. They are also known as stable Paretian laws. We have also mentioned the nature of scaling when the family is discrete. Alternate stable laws include max-stable, min-stable, multiplication-stable, geometric sum-stable, geometric max-stable, geometric min-stable and geometric multiplication-stable, see Mittnik and Rachev (1993).

One may consider the counter parts of sum-stable, max/min-stable and multiplication-stable laws with random sample size. Of these all except the random multiplication-stable laws are already discussed here.

## 10. Applications

Mittnik and Rachev (1993) discuss applications of alternate stable laws in asset returns and their comparison in the context. They also showed that the geometric sum-stable laws dominate all other models. Other applications of this model are discussed in Kozubowski (1999), Kozubowski and Rachev (1994), Kozubowski and Panorska (1999), Kozubowski and Podgorski (2001) and Solomon and Richmond (2002). More applications of alternate stable laws are available in Rachev (1993), Gnedenko and Korolev (1996), Focardi (2002), Makowski (2001), Denuit *et al.* (2000), Kaufman (2001), Kotz *et al.* (2001) and Sajikumar (2002).

Continuous time random walks were introduced in Montroll and Weiss (1965) and are now used in Physics to model anomalous diffusion and relaxation, see Klafter *et al.* (1987), Scher and Lax (1973), Uchaikin and Zolotarev (1999) and Kotuluski and Weron (1996). Kozubowski *et al.* (2002) discusses limit theorems for continuous random walks and shows that the limit is either operator ν-stable or belongs to the generalized domain of normal attraction of an operator stable law. Here the random sample size corresponds to the number of jumps by time $t \geq 0$.

**Remark 6** (see remark.4) By virtue of our Theorem.8 based on Lemmata.2 & 3 we have a class of distributions modeling this number of jumps and this is a flexibility achieved by our operator φ-stable laws over the operator ν-stable laws.

## 11. Free Probability Theory

To the classical probability theory of sums of r.vs there is a parallel in free probability theory or free r.vs; see Voiculescu (1986), Maassen (1992) and Bercovici and Voiculescu (1993). A more recent work in this area is Bercovici and Pata (1999) (with an appendix by Biane) discussing weak convergence, stable laws, domain of attraction and unimodality in the free theory.

*Acknowledgement.* The author wishes to thank Professor Arup Bose for the encouragement and Dr A M Iksanov for the valuable suggestions for the improvement in section.2 and supplying many more references on self-decomposability.

## References

Alice Thomas and Jose, K. K. (2003*a*). Marshall-Olkin Weibull processes, accepted in the *Calcutta Statist. Assoc. Bulletin*.

Alice Thomas and Jose, K. K. (2003*b*). Marshall-Olkin Pareto processes, *Far-East J. Theor. Statist.*, 9, 117 – 132.

Anil, V. (2001). A generalized Poisson distribution and its applications, *J. Kerala, Statist. Assoc.*, 12, 11 – 22.

Arnold, B.C; Robertson, C.A; and Yeh, H.C. (1986): Some properties of a Pareto type distribution, *Sankhya-A*, 404-408.

Balakrishna, N and Jayakumar, K. (1996). Bivariate autoregressive minification process, *J. Appl. Statist. Sci.*, 5, 129 – 141.

Balakrishna, N. and Jayakumar, K. (1997): Bivariate semi-Pareto distributions and processes, *Statist. Papers*, 38, 149 – 165.

Balkema, A. A and Resnick, S. (1977). Max-Infinite divisibility, *J. Appl. Probab*, 14, 309–319.

Barndorff-Nielson, O. E. and Shephard, N. (2000). Modelling by Levy processes for Financial Econometrics, *MPS-RR 2000-16 at www.maphysto.dk*.

Bercovici, H. and Pata, V. (1999) (with an appendix by Biane, P). Stable laws and domains of attraction in free probability theory, *Ann. Math.*, 149, 1023 – 1060.

Bercovici, H. and Voiculescu, D. (1993). Free convolution of measures with unbounded support, *Indiana Univ. Math. J.*, 42, 733 – 773.

Bondesson, L. (1981). Classes of infinitely divisible distributions and densities, *Z. Wahr. Verw. Geb.*, 57, 39 – 71.

Bondesson, L. (1990): Generalized gamma convolutions and complete monotonicity, *Prob. Theor. Rel. Fields*, 85, 181 – 194.

Bose, A; Dasgupta, A. and Rubin, H. (2002). A contemporary review and bibliography of infinitely divisible distributions and processes, *Sankhya-A*, 64, 763 – 819.

Bunge, J. (1996). Composition semi groups and random stability, *Ann. Probab.*, 24, 3, 1476 – 1489.

Cai, I and Kalashnikov, V (2000). NWU property of a class of random sums, *J. Appl. Probab.*, 37, 283 – 289.

Christoph, G. and Schreiber, K. (1998): Discrete Stable random variables, *Statist. Prob. Letters*, 37, 243 – 247.

Denuit, M, Geneset, C and Marceau, E (2002). Criteria for the stochastic ordering of random sums, with Acturial applications, *Scandinavian Journal of Statistics*, 3-16.

Devroye, L. (2001). Algorithms for generating discrete random variables with a given generating function or a given moment sequence, Pre-print (From the internet).

Diedhiou, A.(1998). On the self-decomposability of the half-Cauchy distribution, *J. Math. Anal. Appl.*, 220, 42 – 64.

Dimitrov, B. and Kolev, N. (2002). Beta transformation, Beta type self-decomposability and related characterizations, Preprint.

Dorea, C. C. Y. and Goncalves, C. R. (1999). Asymptotic distribution of extremes of randomly indexed random variables, *Extremes*, 2:1, 95 –109.

Erdogan, M.B. and Ostrovskii, I.V. (1998). Analytic and asymptoptic properties of generalized Linnik probability densities, *J. Math. Anal. Appl.*, 217, 555-578.

Feller, W. (1968). *An Introduction to Probability Theory and Its Applications*, Vol.1, 3rd Edition, John Wiley and Sons, New York.

Feller, W. (1971). *An Introduction to Probability Theory and Its Applications*, Vol. 2, 2nd Edn, John Wiley and Sons, New York.

Focardi, S. M. (2001). Fat tails, scaling and stable laws: A critical look at modeling extremal events in economic and financial phenomena, Discussion paper 2001-02.

Gaver, D.P and Lewis, P.A.W. (1980). First-order auto regressive gamma sequences and point processes, *Adv. Appl. Prob.*, 12, 727 – 745.

Gnedenko, B. V. (1982): On limit theorems for a random number of random variables, *Lecture Notes in Mathematics*, Vol.1021, Fourth USSR-Japan Symp. Proc., Springer, Berlin, 167-176.

Gnedenko, B.V. and Korelev, V.Yu. (1996). *Random Summation, Limit Theorems and Applications*, CRC Press, Boca Raton.

Goldie, C.M. (1967): A class of infinitely divisible distributions, *Proc. Camb. Phil. Soc.*, 63, 1141 – 1143.

Grandall, J. (1991). *Aspects of Risk Theory*, Springer-Verlag, New York.

Grandall, J. (1997). *Mixed Poisson Processes*, Chapman & Hall, London.

Hu, C. Y; Iksanov, A. M; Lin, G.D. and Zakusylo, O. K. (2002). The Hurvitz zeta distribution, Preprint, 29 March 2002.

Iksanov, A. M. and Jurek, Z. J. (2003). Shot noise distributions and self-decomposability, *Stoch. Anal. Appl.*, 21, 593 – 609.

Ismail, M. E. H. and Kelker, D. H. (1979). Special functions, Steiltjes transforms and infinite divisibility, *SIAM J Math. Anal.*, 10, 884 –901.

Jayakumar, K. (1995). The stationary solution of a first order integer valued autoregressive process, *STATISTICA, LV*, 221 – 228.

Jayakumar, K. and Pillai, R. N. (1993). First order autoregressive Mittag-Leffler process, *J. Appl. Probab.*, 30, 462 –466.

Jeanblanc, M; Pitman, J. and Yor, M. (2001). Self-similar processes with independent increments associated with Levy and Bessel processes, *Tech. Report No.608*, Department of Statistics, University of California, Berkeley.

Johnson, N. L.; Kotz, S. and Kemp, A. W. (1992). *Univariate Discrete Distributions*, 2$^{nd}$ edition, Wiley, New York.

Jose, K. K. and Pillai, R. N. (1996). Generalized autoregressive time series modeling in Mittag – Leffler variables, in *Recent Advances in Statistics* (Selected Papers of Prof. R N Pillai), 96 – 103, *(Also presented at the seminar on Recent Trends in Probability and Statistics held in honour of Prof. R N Pillai on 28 June 1993 at University of Kerala, Trivandrum)*.

Jose, K. K. and Alice Thomas (2002). Multivariate minification processes, *STARS Int. J.*, 1, 1-9.

Jose, K. K. and Alice Thomas (2003). Bivariate semi-Pareto processes, accepted in *Metrika*.

Jose, K. K. and Seetha Lakshmi, V. (1999). On geometric exponential distribution and its applications, *J. Ind. Statist. Assoc.*, 37, 51 –58.

Jurek, Z.J. (1997). Self-decomposability: an exception or a rule?, *Ann. Univ. Mariae Curie-Sklodowska, Lubin-Polonia*, L1.1, Section.A, 93 – 106.

Jurek, Z.J. (2000). A note on gamma random variables and Dirichlet series, *Statist. Probab. Lett.*, 49, 387 -392.

Jurek, Z.J. (2001). Remarks on self-decomposability and new examples, *Demonstratio Math.*, XXXIV, 241 - 250.

Jurek, Z.J. (2003). Generalized Levy stochastic areas and self-decomposability, Preprint, April 14.

Kalashnikov, V. (1997*). Geometric Sums : Bounds for Rare Events with Applications*, Klewer Academic Publications, Dordrecht.

Kallenberg, O. (1975). Limits of compound and thinned point processes, *J. Appl. Probab.*, 12, 269 – 278.

Katti, S. K. (1967). Infinite divisibility of integer valued random variables, *Ann. Math. Statist.*, 38, 1306 – 1308.

Kaufman, E. (2001). *Statistical Analysis of Extrme Values, From Insurance, Finance, Hydrology and Other Fields*, Extended 2$^{nd}$ Edition, Birkhauser.

Keilson, J. and Steutel, F.W. (1974). Mixtures of distributions, moment inequalities and measures of exponentiality and normality, *Ann. Probab.*, 2, 112 – 130.

Kemp, A. W. (1978). Cluster size probabilities for generalized Poisson distributions, *Commun. Statist. – Theor. Meth.*, A7 (15), 1433 – 1438.

Klafter, J.; Blumen, A. and Shlesinger, M. F. (1987). Stochastic pathways to anomalous diffusion, *Phys. Rev. A*, 35, 3081 – 3085.

Klebanov, L. B. and Rachev, S. T. (1996). Sums of a random number of random variables and their approximations with ν-accompanying infinitely divisible laws, *Serdica Math. J.*, 22, 471 – 496.

Klebanov, L.B; Maniya, G.M. and Melamed, I.A. (1984). A problem of Zolotarev and analogues of infinitely divisible and stable distributions in the scheme of summing a random number of random variables, *Theor. Probab. Appl.*, 29, 791 – 794.

Knight, F. B. (2001). On the path of an inert object impinged on one side by a Brownian particle, *Probab. Theory Relat. Fields*, 121, 577 – 598.

Kotulski, M. and Weron, K. (1996). Random walk approach to relaxation in disordered systems, *Proceedings of Athens Conference on Applied Probability and Time Series Analysis, Vol. 1*, 379 – 388; *Lecture Notes in Statistics, 114*, Springer, New York .

Kotz, S; Kozubowski, T. J; and Podgorski, K. (2001). *The Laplace Distribution and Generalizations: A Revisit with Applications to Communications, Economics, Engineering, and Finance*, Birkhauser, Boston.

Kozubowski, T. J. (1999). Geometric stable laws: Estimation and applications, *Math. Comput. Modell.*, 29, 241 – 253.

Kozubowski, T. J. (2000). Computer simulation of geometric stable distributions, *J. Computational Appl. Math.*, 116, 221 – 229.

Kozubowski, T. J; Meerschacrt, M. M and Scheffler, H. P. (2002). The operator ν-stable laws, preprint.

Kozubowski, T. J. and Panorska, A. K. (1999). Simulation of geometric stable and other limiting multivariate distributions arising in random summation scheme, *Math. Comput. Modelling*, 29, 255 – 262.

Kozubowski, T,J. and Panorska, A.K. (1996). On moments and tail behavior of ν-stable random variables, *Statist. Prob. Letters*, 29, 307 - 315.

Kozubowski, T,J. and Panorska, A.K. (1998). Weak limits for multivariate random sums, *J. Multi. Anal.*, 67, 398-413.

Kozubowski, T. J; and Podgorski, K. (2001). Asymmetric Laplace laws and modeling financial data, *Math. Comput. Modell.*, 34, 1003 – 1021.

Kozubowski, T.J; and Rachev, S.T. (1994): The theory of geometric stable distributions and its use in modelling financial data, *Europian J. Oper. Res.*, 74, 310 – 324.

Kozubowski, T.J. and Rachev, S.T. (1999a): Univariate geometric stable laws, *J. Computational Anal. Appl.*, 1, 2, 177 – 217.

Kozubowski, T.J. and Rachev, S.T. (1999b): Multivariate geometric stable laws, *J. Computational Anal. Appl.*, 1, 4, 349 – 385.

Lukacs, E. (1970): *Characteristic Functions*, 2[nd] Edition, Griffin, London.

Maassen, H. (1992). Addition of freely independent random variables, *J. Funct. Anal.*, 106, 409 –438.

Makowski A M (2001). On a random sum formula for the busy period of the M/G/∞ queue with applications, Technical Research Report-CSHCN TR 2001-4, *University of Maryland*.

Marshall, A.W; and Olkin, I. (1997): On adding a parameter to a distribution with special reference to exponential and Weibull models, *Biometrika*, 84, 641-652.

McFadden, J. A. (1965). The mixed Poisson process, *Sankhya-A*, 27, 83 – 92.

Meerschacrt, M. M. and Scheffler, H. P. (2001). *Limit Theorems for Sums of Independent Random Vectors*, Wiley New York.

Mittnick, S and Rachev, S T (1993). Modeling asset returns with alternative stable distributions, *Econometric Review*, 12, 261 – 330.

Mohan, N. R. (1998). On geometrically max infinitely divisible laws, *J. Ind. Statist. Assoc.*, 36, 1 –12.

Montroll, E. W. and Weiss, G. H. (1965). Random walks on lattices II, *J. Math. Phys.*, 6, 167 – 181.

Ospina, A. V. and Gerber, H. U. (1987). A simple proof of Feller's characterization of the compound Poisson distributions, *Insurance: Mathematics and Economics*, 6, 63 – 64.

Pillai, R. N. (1990). Harmonic mixtures and geometric infinite divisibility, *J. Ind. Statist. Assoc.*, 28, 87–98.

Pillai, R.N. (1991): Semi-Pareto processes, *J. Appl. Prob.* 28, 461-465.

Pillai, R.N. and Anil, V. (1996): Symmetric stable, α-Laplace, Mittag-Leffler and related laws and processes and the integrated Cauchy functional equation, *J. Ind. Statist. Assoc.*, 34, 97-103.

Pillai, R. N. and Jayakumar, K. (1994) Specialized class-L property and stationary auto regressive process, *Statist. Probab. Lett.*, 19, 51 – 56.

Pillai, R. N. and Jayakumar, K. (1995) Discrete Mittag-Leffler distributions, *Statist. Probab. Lett.*, 23, 271 – 274.

Pillai, R. N; Jose, K.K and Jayakumar, K. (1995). Autoregressive minification processes and the class of distributions of universal geometric minima, *J. Ind. Statist. Assoc.*, 33, 53 – 61.

Pillai, R.N. and Sabu George (1984): A certain class of distributions under normal attraction, *Proc. VI[th] Annual Conf. ISPS*, 107 – 112.

Pillai, R.N. and Sandhya, E. (1990): Distributions with complete monotone derivative and geometric infinite divisibility, *Adv. Appl. Prob..*, 22, 751 – 754.

Pillai, R.N. and Sandhya, E. (1996). Geometric sums and Pareto law in reliability theory, *IAPQR Trans.*,21, 137-142.

Pillai, R. N. and Sandhya, E. (2001). A characteristic property of inverse Gaussian law, *Statistical Methods*, 3, 36 – 39.

Pillai, R. N. and Satheesh, S. (1992). α-inverse Gaussian distributions, *Sankhya-A*, 54, 288 – 290.

Rachev, S T (1993). Rate of convergence for maxima of random arrays with applications to stock returns, *Statistics and Decisions*, 11, 279 – 288.

Rachev, S.T. and Resnick, S. (1991) Max geometric Infinite Divisibility and Stability, *Comm. Statist, -Stoch. Models*, 7, 2, 191 – 218.

Ramachandran, B. (1997): On geometric stable laws, a related property of stable processes, and stable densities of exponent one, *Ann. Inst. Statist. Math.*, 49, 2, 299 – 313.

Riedel, M. (1980*a*). Representation of the characteristic function of a stochastic integral, *J. Appl. Probab.*, 17, 448 – 455.

Riedel, M. (1980*b*). Characterization of stable processes by identically distributed stochastic integrals, *J. Appl. Probab.*, 689 – 709.

Romano, J. P. and Siegel, A. F. (1986). *Counterexamples in Probability and Statistics*, Chapman and Hall, New York.

Sajikumar, V. R. (2002). *Distributions Related to Mittag-Leffler Function and Modeling Financial Data*, Ph.D. Thesis (unpublished), University of Kerala.

Sandhya, E. (1991*a*). *Geometric Infinite Divisibility and Applications*, Ph.D. Thesis (unpublished), University of Kerala, January 1991.

Sandhya, E. (1991*b*). On geometric infinite divisibility, *p*-thinning and Cox processes, *J. Kerala Statist. Assoc.*, 7, 1-10.

Sandhya, E. (1996). On a generalization of geometric infinite divisibility, *Proc. 8[th] Kerala Science Congress*, January 1996, 355 – 357.

Sandhya, E. and Pillai, R.N. (1999). On geometric infinite divisibility, *J. Kerala Statist. Assoc.*,10, 1–7.

Sandhya, E. and Satheesh, S. (1996*a*). On the membership of semi-α-Laplace laws in class-L, *J. Ind. Statist. Assoc.*, 34, 77 – 78.

Sandhya, E. and Satheesh, S. (1996*b*): On distribution functions with completely monotone derivative, *Statist. Prob. Letters*, 27, 127 – 129.

Sandhya, E. and Satheesh, S. (1997): On exponential mixtures, mixed Poisson processes and generalized Weibull and Pareto models, *J. Ind. Statist. Assoc.*, 35, 45 – 50.

Sato, K. (2001). Subordination and self-decomposability, *Statist. Probab. Lett.*, 54, 317 – 324.

Satheesh, S. (2001*a*). *Stability of Random Sums and Extremes*, Ph.D. Thesis (unpublished), Cochin University of Science and Technology, July 2001.

Satheesh, S. (2001*b*). Another Look at Random Infinite Divisibility, submitted in January 2002. (Results were presented at the *International Conference on Statistics in Industry and Business, Cochin, 2-4 of January 2003.*).

Satheesh, S. (2002). Aspects of Randomization in Infinitely Divisible and Max-Infinitely Divisible Laws, *ProbStat Models*, 1, June – 2002, 7 – 16.

Satheesh, S; and Nair, N.U. (2000): On the stability of geometric extremes, Submitted, (Presented at the *International Conference on Order Statistics and Extreme Values: Theory and Applications, University of Mysore, December 2000.*).

Satheesh, S. and Nair, N.U. (2002*a*). Some classes of distributions on the non-negative lattice, *J. Ind. Statist. Assoc.*, 40, 41 – 58.

Satheesh, S. and Nair, N.U. (2002*b*). A Note on maximum and minimum stability of certain distributions, *Calcutta Statist. Assoc. Bulletin*, 53, 249 -252.

Satheesh, S; Nair, N U and Sandhya, E (2002). Stability of random sums, *Stochastic Modeling and Applications*, 5, 17 –26.

Satheesh, S. and Sandhya, E. (1997). Distributions with completely monotone probability sequences, *Far East J. Theor. Statist.*, 1, 69-75.

Satheesh, S. and Sandhya, E. (2002). Infinite divisibility and max-infinite divisibility with random sample size, accepted in *Statistical Methods.* (*Results were presented at the International Conference on Statistics in Industry and Business, Cochin, 2-4 of January 2003*).

Scher, H. and Lax, M. (1973). Stochastic transport in a disordered solid - 1, *Theory. Phys. Rev., B* (7), 4491 – 4502.

Seethalekshmi, V. and Jose, K. K. (2001). A subordinated process with geometric exponential operational time, *Calcutta Statist. Assoc. Bulletin*, 51, 273 – 275.

Seethalekshmi, V. and Jose, K. K. (2002). Geometric Mittag-Leffler tailed autoregressive processes, *Far-East J. Theor. Statist.*, 6, 147 – 153.

Seethalekshmi, V. and Jose, K. K. (2003). Geometric Mittag-Leffler distribution and applications, accepted in *J. Appl. Statist. Science.*

Seethalekshmi, V. and Jose, K. K. (2004). Geometric α-Laplace distribution and auto regressive processes, accepted in *Statistical Papers.*

Seshadri, V. (1993). *The Inverse Gaussian Distribution, A Case Study in Exponential Family of Distributions*, Oxford University Press, New York.

Shaked, M. (1975): On the distributions of the minimum and of the maximum of a random number of i.i.d random variables, *In Statistical Distributions in Scientific Work*, Editors, Patil, G.P; Kotz, S; and Ord, J.K., D.Reidel Publishing Company, Dordrecht, Holland, 363-380.

Shanbhag, D.N. and Sreehari, M. (1977): On certain self-decomposable distributions, *Z. Wahr. Verw. Geb.* , 38, 217 – 222.

Shanbhag, D.N and Sreehari, M. (1979): An extension of Goldie's result and further results on infinite divisibility, *Z. Wahr. Verw. Geb.*, 47,19–25.

Shanbhag, D.N; Pestana, D. and Sreehari, M. (1977). Some further results in infinite divisibility, *Math. Proc. Camb. Phil. Soc.*, 82, 289 – 295.

Silvestrov, D.S. and Teugels, J.L. (1998). Limit theorems for extremes with random sample size, *Adv. Appl. Probab.*, 30, 777 – 806.

Solomon, S. and Richmond, P. (2002). Stable power laws in variable economics; Lotka–Volterra implies Pareto–Zipf, *Euro. Phys. J. – B*, 27, 257 – 261.

Sreehari, M (1995). Maximum Stability and a generalization, *Statist. Probab. Lett.*, 23, 339 – 342.

Steutel, F.W. (1973). Some recent results in infinite divisibility, *Stoch. Proc. Appl.*, 1, 125 – 143.

Steutel, F.W. (1979). Infinite divisibility in theory and practice, *Scand. J. Statist.*, 6, 57 – 64.

Steutel, F.W. (1969): Note on completely monotone densities, *Ann. Math. Statist.*, 40, 1130 – 1131.

Steutel, F.W. and van Harn, K. (1979). Discrete analogues of self-decomposability and stability, *Ann. Prob.*, 7, 893-899.

Sundt, B. (2000). Multivariate compound Poisson distributions and infinite divisibility, *Austin Bulletin*, 30, 305 – 308.

Thorin, O. (1977): On the infinite divisibility of the Pareto distribution, *Scand. Actuarial J.*, 31 – 40.

Uchaikin, V. V. and Zolotarev, V. M. (1999). *Chance and Stability, Stable Distributions and Their Applications*, VSP, Utrecht.

Voiculescu, D. (1986). Addition of certain non-commuting random variables, *J. Funct. Anal.*, 66, 323 – 346.

Voorn, W. J. (1987). Characterizations of the logistic and log-logistic distributions by extreme value related stability with random sample size, *J. Appl. Probab.,* 24, 838 - 851.

Voorn, W. J. (1989). Stability of extremes with random sample size, *J. Appl. Probab.*, 27, 734 – 743.

S Satheesh

Neelolpalam

S. N. Park Road

Trichur – 680 004, **India**.

e-mail: ssatheesh@sancharnet.in

The following papers are also available at http://arxiv.org as:

Satheesh, S. (2001*b*). Another Look at Random Infinite Divisibility; math.PR/0304499.

Satheesh, S. (2002). Aspects of Randomization in Infinitely Divisible and Max-Infinitely Divisible Laws, *ProbStat Models*, 1, June – 2002, 7 – 16; math.PR/0305030.

Satheesh, S; and Nair, N.U. (2000). On the stability of geometric extremes; math.PR/0305060.

Satheesh, S. and Nair, N.U. (2002*b*). A note on maximum and minimum stability of certain distributions; math.PR/0305059.

Satheesh, S. and Sandhya, E. (2002). Infinite divisibility and max-infinite divisibility with random sample size; math.PR/0305045.